\documentclass[10pt]{amsart}
\usepackage{amsmath}
\usepackage{amsthm}
\usepackage{amsfonts}
\usepackage{amssymb}
\usepackage{latexsym}
\usepackage{enumerate}

\newcommand{\beq}{\begin{eqnarray}}
\newcommand{\eeq}{\end{eqnarray}}
\newcommand{\beqno}{\begin{eqnarray*}}
\newcommand{\eeqno}{\end{eqnarray*}}
\newcommand{\bR}{\mathbb R}

\newcommand{\cG}{\mathcal G}

\newcommand{\cF}{\mathcal F}

\newcommand{\cA}{\mathcal A}
\newcommand{\bM}{\bar{M}}

\newcommand{\midint}{- \mskip-18mu \int}
\newcommand{\Hom}{\mbox{Hom}}
\newcommand{\eps}{\varepsilon}
\renewcommand{\phi}{\varphi}
\newcommand{\dx}{\, dx}

\allowdisplaybreaks

\def\mvint_#1{\mathchoice
  {\mathop{\vrule width 6pt height 3 pt depth -2.5pt
        \kern -8pt \intop}\nolimits_{\kern -3pt #1}}%
 {\mathop{\vrule width 5pt height 3 pt depth -2.6pt
                  \kern -6pt \intop}\nolimits_{#1}}%
 {\mathop{\vrule width 5pt height 3 pt depth -2.6pt
                  \kern -6pt \intop}\nolimits_{#1}}%
 {\mathop{\vrule width 5pt height 3 pt depth -2.6pt
                  \kern -6pt \intop}\nolimits_{#1}}}

\newtheorem{theorem}{Theorem}[section]
\newtheorem{proposition}[theorem]{Proposition}
\newtheorem{Lem}[theorem]{Lemma}

\theoremstyle{definition}
\newtheorem{definition}[theorem]{Definition}
\newtheorem*{remark}{Remark}

\begin{document}
\title[Regularity results by $A$-harmonic approximation]{Regularity for minimizers of functionals with nonstandard growth by
A-harmonic approximation}

\author{Jens Habermann}
\address{Jens Habermann, Institute for mathematics\\
Friedrich-Alexander University\\ Bismarckstr. 1 1/2\\ 91054
Erlangen\\ Germany;}
\email{habermann@mi.uni-erlangen.de}

\author{Anna Zatorska--Goldstein}
\address{Anna Zatorska--Goldstein, Institute of Applied Mathematics and Mechanics\\
University of Warsaw\\ Banacha 2\\ 00-913 Warsaw\\ Poland;}
\email{azator@mimuw.edu.pl}
\thanks{A.Z.-G. is partially supported by MEiN grant no 1PO3A 005 29 and by Alexander von
Humboldt Foundation}
\date{\today}

\begin{abstract}
We prove partial regularity for minimizers of quasiconvex  functionals of the
type $\int_\Omega f(x,Du) \,dx$ with $p(x)$ growth with respect to the
second variable. The proof is direct and it uses a method of $\cA$--harmonic
approximation.
\end{abstract}

\maketitle

\section{Introduction}
In this paper we study the regularity properties of local minimizers of
a variational functional
\begin{equation*}
     \cF[u] = \int_{\Omega} f(x,Du)\dx,
\end{equation*}
where $u : \Omega \to \bR^N$, $\Omega$ is a bounded domain in $\bR^n$ and
the integrand $f:\Omega \times \Hom(\bR^n;\bR^N) \to \bR$ satisfies the growth
condition of the type
 \begin{equation*}
    f(x,A)  \approx (1 + |A|^2)^{p(x)/2},
\end{equation*}
for $p:\Omega \to (1, \infty)$ beeing a H\"older continuous function. For the
precise statement of the conditions see section 2.

\medskip

\begin{definition}
A function $u\in W_{loc}^{1,1}(\Omega;\bR^n)$ is called a local minimizer of the functional $\cF$ if
$|Du|^{p(x)} \in L_{loc}^1(\Omega)$ and
\begin{equation*}
     \cF[u] \le \cF[u + \varphi],
\end{equation*}
for all $\varphi \in W^{1,1}_0(\Omega;\bR^N)$ with compact support in $\Omega$.
\end{definition}

The main statement is the following
\begin{theorem}
Let $u \in W^{1,1}_{loc}(\Omega;\bR^N)$ be a local minimizer of the functional
$\cF$ fulfilling the assumptions {\bf A1} -- {\bf A3} (see page \pageref{assumptions}). 
Let $\gamma_2$ be an upper bound for the exponent $p$ and assume that the modulus of
continuity $\omega$ of $p$ satisfies the condition
\begin{equation}\label{cont.strong}
     \omega(\rho) \le L\rho^{\alpha},
\end{equation}
for some $L>1$, $\alpha \in (0,1]$ and all $\rho < 1$. Then there
exists an open subset $\Omega_0 \subset \Omega$ with ${\mathcal
L}^n(\Omega \setminus \Omega_0) = 0$ such that $Du\in
C^{0,\beta}_{loc}(\Omega_0)$ with $\beta = \min \{ 1, \frac{2}{\gamma_2}\}
\frac{\alpha}{4}$.
\end{theorem}

\medskip

The proof of the theorem (with different $\beta$) was done by E.
Acerbi and G. Mingione in 2001 \cite{acerbi-m2}. The key step is to
establish a certain excess--decay estimate for the so called excess
function $\Phi$, which is defined as
\begin{equation} \label{excess}
     \Phi \equiv \Phi(x_0,\rho,A) \equiv \left( \midint_{B_{\rho}(x_0)} |V_{p_2}(Du) -
     V_{p_2} (A)|^2\dx \right)^{1/2},
\end{equation}
where $p_2$ denotes the maximal exponent $p(x)$ in a neighbourhood of $x_0$ and
with $V_p : \bR^k \to \bR^k$ given by
\begin{equation}\label{V}
V_p(\xi) = \left(1+|\xi|^2\right)^{(p-2)/4} \xi.
\end{equation}
The function $\Phi$ provides an integral measure of the oscillations of the
gradient $Du$ in a ball $B_\rho$.  The excess-decay estimate leads to H\"older
continuity of $Du$ in $B_\rho$ via the integral characterization of H\"older
continuous functions due to Campanato (see \cite{campanato}).
The excess-decay estimate was established by Acerbi and Mingione in an indirect
way, using the blow--up technique.

We present here a more direct proof of the result. Applying the
variational principle of Ekeland we obtain a comparison function,
i.e. a function which is an almost minimizer of the functional with
frozen $x$--coefficient and which is close to our local minimizer in
an appropriate Sobolev norm. Having such a comparison function at
hand, we are able to use the results for almost minimizers of
variational functionals with constant $p$ growth. In particular, we
obtain straightforward a Cacciopoli type inequality for local
minimizers of the functional $\cF$. Then, instead of blow--up
arguments we use a method of $\cA$--harmonic approximation to obtain
an excess-decay estimate.

The method originates in a work of L. Simon. It is based on the fact
that one is able to obtain a good approximation of a function $w \in
W^{1,2}(B;\bR^N)$, which is approximately $\cA$--harmonic in a
certain sense by an $\cA$--harmonic function $h\in W^{1,2}(B;\bR^N)$,
in both the $L^2$--topology and the weak topology in $W^{1,2}$. Here $h$
is called $\cA$--harmonic on $B$ if there holds
$$
\int_B \cA(Dh, D\phi) \dx = 0 \quad \text{for any $\phi \in C_0^1(B;\bR^N)$},
$$
where $\cA$ is a bilinear form on $Hom(\bR^n; \bR^N)$ which is
(strongly) elliptic in the sense of Legendre--Hadamard, i.e. for all
$\eta \in \bR^n$ and $\xi \in \bR^N$ there holds:
$$
\cA(\eta \otimes \xi, \eta \otimes \xi) \geq \kappa |\eta|^2 |\xi|^2.
$$

We can assume that our exponent function $p$ varies in some bounded
interval $[\gamma_1, \gamma_2]$. However we have to be able to
consider both cases $\gamma_2 \geq 2$ and $1< \gamma_2 <2$. For $p
\geq 2$ it is straightforward to adapt the standard $\cA$--harmonic
approximation lemma by using the $L^2$--theory combined with a standard
Sobolev inequality. For $1<p<2$ we do not have the access to the
$L^2$--theory for functions in $W^{1,p}$ but still it is possible to
generalize the approximation lemma directly. This was done in
\cite{duzaar-g-k}.

Apart from the fact that the proof is very clear, the method provides
better control of the constants. We have to admit that the proof of the
$\cA$--approximation lemma itself is done by contradiction. In
direct way we show that a local minimizer is approximately
$\cA$--harmonic, where $\cA = D^2f(x_0, (Du)_{x_0,\rho})$. The
$\cA$--harmonic approximation lemma guarantees the existence of a
certain constant which is, admittedly not in an explicit form,
determined by a property of constant coefficient elliptic systems,
and we will later use the constant in the regularity proof. This
constant, however, does not have an influence on the final H\"older
exponent of $Du$. We should mention that by our method we end up with
a final H\"older exponent
\begin{equation*}
     \beta \equiv \min\{1,2/\gamma_2\} \frac{\alpha}{4}
\end{equation*}
for the function $Du$,
where $\alpha$ denotes the H\"older exponent of the exponent 
function $p$ and $\gamma_2$ is the global bound for $p$. This is, 
in fact, a slightely better result than stated in \cite{acerbi-m2}.

\bigskip

In order to obtain  regularity results for local minimizers of the functional
$\cF$ we have to assume some continuity properties of the exponent $p$.  The minimal condition about the
modulus of continuity of $p$ is that 
\begin{equation}\label{cont.weak}
     \limsup_{\rho \to 0} \omega(\rho) \log\left(\frac{1}{\rho}\right) <
     \infty.
\end{equation}
Dropping this assumption in general causes the loss of any type of regularity of
minimizers (see \cite{zhikov2}).
By the result of Zhikov \cite{zhikov}, \eqref{cont.weak} is sufficient to obtain higher integrability of
the gradient of a minimizer. However, it is not sufficient to obtain further
regularity.  Acerbi and Mingione \cite{acerbi-m2} proved $C^{0,\alpha}$
regularity for minimizers for every $\alpha <1$, provided the modulus of
continuity satisfies an assumption
\begin{equation}\label{cont.stronger}
     \limsup_{\rho \to 0} \omega(\rho) \log\left(\frac{1}{\rho}\right) =0,
\end{equation}
which is in accordance with the theory of functionals with constant $p$--growth
where an additional continuity assumption with respect to $x$ is required to
reach any exponent $\alpha <1$. 
In order to prove $C^{1,\alpha}$ regularity of minimizers, in constant $p$ case
(both for $p\in (1,2)$ and $p \geq 2$) assumption \eqref{cont.stronger} is not
sufficient. In fact one needs either that the modulus of continuity satisfies the so called Dini condition or that $p$ itself is H\"older continuous function, i.e.
\begin{equation*}
     \omega(\rho) \approx \rho^{\alpha}.
\end{equation*}
This condition was assumed by Acerbi and Mingione in the original proof of the
result of H\"older continuity of the gradient of a minimizer in $p(x)$ case.

\section{Setting}

We impose the following structure conditions for the functional $\cF$:

\noindent{\bf A1 (growth):}\label{assumptions} the function $f(x,\cdot)$ is of the class $C^2$ and there exist constants
$\mu \in (0,1]$, $L \ge 1$ such that for all $x \in \Omega$ and  $A \in \Hom(\bR^n;\bR^N)$
 we have
\begin{equation}\label{growth}
     L^{-1} |A|^{p(x)} \le f(x,A) \le L(\mu^2 +
     |A|^2)^{p(x)/2},
\end{equation}
where $p:\Omega \to (1,\infty)$ is a continuous function;

\medskip

\noindent{\bf A2 (quasi-convexity):} the function $f(x,\cdot)$ is (strictly) quasi--convex i.e.
\begin{multline}\label{ellipt.}
     \int_{B_\rho(x_0)} (f(x_0,A+ D\varphi) - f(x_0,A))\dx \\
     \ge \frac{1}{L} \int_{B_\rho(x_0)} (\mu^2 +
     |A|^2 + |D\varphi|^2)^{\frac{p(x_{0})-2}{2}} |D\varphi|^2\dx,
\end{multline}
for all $x_0 \in \Omega$, $B_\rho(x_0) \Subset \Omega$,
$A \in \Hom(\bR^n;\bR^N)$ and $\varphi \in C_0^{\infty}(B_\rho(x_0);\bR^N)$;

\medskip

\noindent{\bf A3 (continuity):} the function $f$ satisfies the following continuity condition with respect to
the first variable
\begin{multline} \label{continuity x}
     |f(x,A) - f(x_0,A)| \le L\omega(|x-x_0|)
    \left[ (\mu^2+|A|^2)^{p(x)/2} +
     (\mu^2+|A|^2)^{p(x_0)/2}\right] \cdot \\
     \cdot \Big[1+|\log(\mu^2+|A|^2)|\Big],
\end{multline}
for all $x,x_0 \in \Omega$ and $A \in \Hom(\bR^n,\bR^N)$, where 
$\omega: (0,\infty) \to (0,\infty)$ is a modulus of continuity for the function
$p$, i.e. a non-decreasing, continuous
function with $ \lim_{R \to 0} \omega(R) = 0$ and
\begin{equation*}
     |p(x) - p(y)| \le \omega(|x-y|) \quad \mbox{ for all}\ x,y \in \Omega.
\end{equation*}
Remark that \eqref{continuity x} is a natural  condition for $f(x,A) = (1+|A|^2)^{p(x)/2}$. 

\bigskip

Since $f(x,\cdot)$ is quasi-convex and satisfies the growth
condition \eqref{growth} it is well known that there exists a
constant $c=c(n,N,p(\cdot ),L)$ such that the first derivatives of $f$
satisfy the growth condition
\begin{equation}\label{growth_Df}
|Df(x_0,A)| \le c(1+|A|^{p(x_0)-1}).
\end{equation}
We do not assume an explicit growth condition for the second
derivatives of $f$. For our purposes it is sufficient that for any
$M>0$ there exists a constant $K_{M,x_0}>0$ such that for $A\in
\Hom(\bR^n; \bR^M)$
\begin{equation}\label{growth_D2f}
\sup_{|A|\le M} |D^2f(x_0,A)| \le K(M,x_0).
\end{equation}
Condition \eqref{growth} implies also the existence of a modulus
of continuity of $D^2f(x,\cdot)$ on compact subsets of
$\Hom(\bR^n;\bR^N)$ i.e. for any given constant $M>0$
\begin{equation}\label{cont_D2f}
     |D^2f(x_0,A) - D^2f(x_0,B)| \le \nu_{M,x_0}(|A-B|),
\end{equation}
for any $A,B \in \Hom(\bR^n,\bR^N)$ with $|A|, |B| \le M+1$.

\begin{remark}
Since our results are of the local nature we will  assume that there exist $1 <
\gamma_1 \le \gamma_2 < \infty$ such that
\begin{equation*}
     \gamma_1 \le p(x) \le \gamma_2  \quad \mbox{ for all}\ x \in \Omega
\end{equation*}
It follows that $K(M,x)$ in \eqref{growth_D2f} may be chosen independently
of $x$. We will therefore omit its dependence on $x$ and write $K(M)$.
\end{remark}

\medskip

\noindent
{\bf Some notation:} Within the whole paper we will write $B(x_0,\rho)$ for the open ball
with centre $x_0 \in \bR^n$ and radius $\rho$. Furthermore we write
\begin{equation*}
     (u)_{x_0,\rho} \equiv \midint_{B(x_0,\rho)} u\, dx \equiv \frac{1}{|B(x_0,\rho)|}
     \int_{B(x_0,\rho)} u\, dx
\end{equation*}
for the mean value of the function $u$ on the ball $B(x_0,\rho)$.
>From time to time we just write $B_{\rho}(x_0)$, or if the center is
clear from the context, $B_{\rho}$ instead of $B(x_0,\rho)$. The
same we do with the notation for the mean value, i.e. we just write
$(u)_{\rho}$ instead of $(u)_{x_0,\rho}$. Concerning the constants
appearing in the proofs we remark that they may change from line to
line. If a constant will be important for the proceeding of the
proofs, we will indicate this in an obvious way. From time to time
for clearness we will not show the dependencies of the constants
within the estimates, but at the end of them.

\section{Basic tools}

\subsection{Higher integrability}
We start with a higher integrability result due to Zhikov, which in the
form of the following statement appears in \cite{acerbi-m1}.
\begin{Lem}\label{High.Int.u}
Let $u \in W^{1,1}_{loc}(\Omega,\bR^N)$ with $|Du|^{p(\cdot )} \in
L^1_{loc}$ be a local minimizer of the functional
\begin{equation*}
     w \mapsto \int_\Omega f(x,Dw(x))\dx,
\end{equation*}
where $f$ satisfies the growth and ellipticity conditions
\eqref{growth}, \eqref{ellipt.} and assumption \eqref{cont.weak} on the modulus
of continuity $\omega$ holds. Furthermore assume that
\begin{equation*}
     \int_\Omega |Du|^{p(x)}\dx \le M < \infty.
\end{equation*}
Then there exist an exponent $\delta =
\delta(n,\gamma_1,\gamma_2,L,M) > 0$, a constant $c = \linebreak
c(n,\gamma_1,\gamma_2,L,M)$ and a radius $R_0 = R_0(n,
\gamma_1,\omega(\cdot ))$ such that for any radius $R \le R_0$ there holds
\begin{equation*}
     \left(\midint_{B_{R/2}} |Du|^{p(x)(1+\delta)}\dx\right)^{\frac{1}{1+\delta}}
     \le c\left(\midint_{B_R} |Du|^{p(x)}\dx + 1\right).
\end{equation*}
\end{Lem}

The next lemma is an up-to-the-boundary result.  The version we present here,
with balls of the same size on both sides of the inequality was proved in
\cite{jens} in higher order case. The original proof (with a ball of double
radius on the right hand side) can be found in
\cite{acerbi-f}, \cite{acerbi-m2}, \cite{carozza-f-m}.

\begin{Lem}[Higher integrability up to the boundary]\label{HI.bdry}
Let $B(x_0,\rho) \Subset \Omega$ and $p$ be a constant such that $1\le\gamma_1 \le p \le
\gamma_2$. Assume $g:\Omega \times
\bR^{nN} \to \bR$ is continuous and for all $z \in \bR^{nN}$ there holds
\begin{equation}\label{growth.g}
     L^{-1}|z|^p \le g(x,z) \le L\left(|z|^p + a(x)\right),
\end{equation}
with $L \ge 1$, $0 < a \in L^{\gamma}
\left(B_{\rho}\right)$, $\gamma >
1$.

Let $h \in W^{1,q}\left(B_{\rho}\right)$ with $q>p$ and $v$ be a
solution of the Dirichlet problem
\begin{equation}
     \min\left\{\int_{B_{\rho}} g\left( x,Dw\right)\dx,\ \ \ w \in
     h + W_0^{1,p}\left(B_{\rho}\right) \right\}.
\end{equation}
Then there exists $\eps = \eps\left(\gamma_1, \gamma_2, L,m\right)
\in (0,m)$ with $m = \min\left\{\gamma - 1, \frac{q}{p}-1\right\}$
and a constant $c \equiv c\left(\gamma_1, \gamma_2,L\right)$ such
that
\begin{multline*}
     \left( \midint_{B_{\rho}}|Dv|^{p(1+\eps)} \dx \right)^{\frac{1}{p(1+\eps)}}
     \le c \Biggr[ \left(\midint_{B_{\rho}} |Dv|^p \dx
     \right)^{\frac{1}{p}} \\
     + \left( \midint_{B_{\rho}} |Dh|^{p(1+m)}\dx \right)^{\frac{1}{p(1+m)}}
     + \left( \midint_{B_{\rho}}a^{1+m}\dx \right)^{\frac{1}{p(1+m)}} \Biggr] .
\end{multline*}
\end{Lem}

\subsection{Ekeland variational principle}

In order to obtain a comparison function, i.e. an almost minimizer
of the functional with frozen coefficients we apply a well known
variational principle of Ekeland (see \cite{ekeland}).

\begin{Lem} \label{eke_orig}
Let $(X,d)$ be a complete metric space and $\cG
: X \to (-\infty,+\infty]$ a lower
semicontinuous functional such that $\inf_X
\cG$ is finite. Given $\eps>0$ let $u\in X$ be
such that $\cG(u) \le \inf_X \cG +\eps$. Then there
exists $w \in X$ such that
\begin{eqnarray*}
d(w,u) &\leq& 1, \\
\cG(w) &\leq& \cG(u),\\
\cG(w) &\leq& \cG(v) +\eps d(v,w), \quad
\text{for any $v\in X$.}
\end{eqnarray*}
\end{Lem}

\subsection{Algebraic properties of the function $V_p$}

Let the function $V \equiv V_p:\bR^k \to \bR^k$ be defined by
\begin{equation}\label{Def.V_p}
     V_p(z) = \left(1+|z|^2\right)^{\frac{p-2}{4}}z.
\end{equation}
We recall algebraic properties of the function
$V_p$ (for a proof of the properties see
e.g. \cite{carozza-f-m}).

\begin{Lem}\label{prop.V}
Let $p>1$ and let $V \equiv V_p: \bR^k \to \bR^k$
be as in \eqref{Def.V_p}. Then for any $z, \eta\in \bR^k$ there holds

\begin{enumerate}[i)]
\item $|V(tz)| \leq \max\{t,t^{p/2}\}|V(z)|$,\quad for any  $t>0$;
\item $|V(z+\eta)| \leq c\Big(|V(z)| +|V(\eta)|\Big)$;
\item
$$
c^{-1}|z-\eta| \leq
\frac{|V(z)-V(\eta)|}{(1+|z|^2+|\eta|^2)^{(p-2)/4}}
\leq c|z-\eta|;
$$
\noindent Moreover for any $z\in \bR^k$
\item
\begin{align*}
\text{if $p\in(1, 2)$:} &\quad \frac{1}{\sqrt{2}}\min\{|z|,|z|^{p/2}\} \leq |V(z)| \leq
\min\{|z|,|z|^{p/2}\}; \\
\text{if $p\geq 2$:} &\quad \max\{|z|,|z|^{p/2}\} \leq |V(z)| \leq
\sqrt{2}\max\{|z|,|z|^{p/2}\};
\end{align*}
\item
\begin{align*}
\text{if  $p \in(1,2)$:}& \quad |V(z)-V(\eta)| \leq c |V(z-\eta)|, \quad
\text{for any $\eta \in \bR^k$}; \\
\text{if  $p \geq 2$:}& \quad |V(z)-V(\eta)| \leq c(M) |V(z-\eta)|, \quad
\text{for $|\eta| \leq M$};
\end{align*}
\item
\begin{align*}
\text{if  $p \in(1,2)$:}& \quad |V(z-\eta)| \leq c(M)|V(z)-V(\eta)|, \quad
\text{for $|\eta| \leq M$}; \\
\text{if  $p \geq 2$:}& \quad |V(z-\eta)| \leq c|V(z)-V(\eta)|, \quad
\text{for any $\eta \in \bR^k$};
\end{align*}
\end{enumerate}
with $c(M), c\equiv c(k,p) >0$. If $1< \gamma_1 \leq p \leq
\gamma_2$ all the constants $c(k,p)$ may be replaced by a single constant
$c\equiv c(k,\gamma_1, \gamma_2)$.
\end{Lem}

\subsection{$\cA$--harmonic approximation and a priori estimates for
$\cA$-harmonic functions}

The key ingredient of the proof is the  following $\cA$--harmonic approximation
lemma.  The proof for the case $p \geq 2$  can be found in \cite{duzaar-s}. The
case $ 1<p<2$ has been proved in \cite{duzaar-g-k}.

\begin{Lem}\label{A-approx}
Let $p>1$ and $\kappa$, $K$ be positive constants. Then for any $\eps>0$ there
exists $\delta=\delta(n,N,\kappa,K,\eps)\in(0,1]$ with the following
property: for any bilinear form $\cA$ on $\Hom(\bR^n;\bR^N)$ which
is elliptic in the sense of Legendre--Hadamard with ellipticity
constant $\kappa$ and upper bound $K$ and  for any $v\in
W^{1,p}(B_\rho(x_0); \bR^N)$ satisfying
\begin{equation*}
\mvint_{B_\rho(x_0)} |V_p(Dv)|^2 dx \le \gamma^2 \le 1 \quad
\text{and}
\end{equation*}
\begin{equation*}
\mvint_{B_\rho(x_0)} \cA(Dv,D\phi) dx \le \gamma \delta
\sup_{B_\rho(x_0)} |D\phi| \qquad \text{ for all $\phi\in
C_0^1(B_\rho(x_0);\bR^N)$},
\end{equation*}
there exists an $\cA$-harmonic function $h$ satisfying
\begin{equation*}
\mvint_{B_\rho(x_0)} |V_p(Dh)|^2 dx \le 1 \quad \text{and} \quad
\mvint_{B_\rho(x_0)} \left|V_p\left(\frac{v-\gamma
h}{\rho}\right)\right|^2 dx \le \gamma^2 \eps.
\end{equation*}
\end{Lem}

\medskip

In Section \ref{excess-dec} we will use  a priori estimates for solutions of linear elliptic systems
of second order with constants coefficients (see e.g. \cite{carozza-f-m} and
\cite{duzaar-g-k}).

\begin{Lem}\label{A-harm}
Let $h \in W^{1,1}(B_\rho(x_0);\bR^N)$ be an ${\mathcal A}$-harmonic function, i.e.
\begin{equation*}
     \int_{B_\rho(x_0)} {\mathcal A}(Dh,D\varphi)\, dx = 0,
\end{equation*}
for any $\varphi \in C_0^1(B_{\rho}(x_0);\bR^N)$, where ${\mathcal A}
\in \Hom(\bR^n; \bR^N)$ is elliptic in the sense of
Legendre-Hadamard with ellipticity constant $\kappa$ and upper bound
$K$. Then $h \in C^{\infty}(B_\rho(x_0);\bR^N)$ and
\begin{equation*}
     \rho \sup_{B_{\rho/2}(x_0)} |D^2 h| + \sup_{B_{\rho/2}(x_0)} |Dh| \le
     c_a \mvint_{B_\rho(x_0)} |Dh|\, dx,
\end{equation*}
where the constant $c_a$ depends only on $n,N,\kappa$ and $K$.
\end{Lem}

\section{Preliminary results}
As remarked before, since our results are of the local nature we will  assume that there exist $1 <
\gamma_1 \le \gamma_2 < \infty$ such that
\begin{equation*}
     \gamma_1 \le p(x) \le \gamma_2  \quad \mbox{ for all}\ x \in \Omega,
\end{equation*}
and moreover
\begin{equation*}
     \int_\Omega |Du|^{p(x)}\dx < \infty.
\end{equation*}
Let $\delta$ be the higher integrability exponent from Lemma
\ref{High.Int.u} and let from now on
the radius $R$ be so small that
\begin{equation*}
     \omega(R) \le \frac{\delta}{4}.
\end{equation*}
Subsequently we will always assume that $\rho \le R$. Take a ball $B_{2\rho}(x_0)$ and define
\begin{equation*}
     p_1 := \inf\{p(x): x \in B_{2\rho}(x_0)\},\quad p_2 := \sup\{p(x):x \in
     B_{2\rho}(x_0)\}.
\end{equation*}
Let furthermore $x_m \in \overline{B_{2\rho}(x_0)}$ be the point,
where the function $p$ reaches the value $p_2$, i.e. $p_2 \equiv
p(x_m)$. Then by $p_2 - p_1 \le \omega(R) \le \delta/4$ we get
\begin{equation}\label{Est.p1p2}
     p_2(1+\delta/4) \le p_1(1+\delta) \le p(x)(1+\delta).
\end{equation}

\subsection{Comparison}

\begin{proposition} \label{ekeland}
Let $B_{2\rho}(x_0) \Subset \Omega$ and
assume that $(|Du|^{p_2})_{x_0,2\rho} \leq \bM < \infty$. Then
there exist
a constant $C(\bM)=C(\bM, \gamma_1, \gamma_2, L, \alpha)>0$
and a function $w\in u + W^{1,p_2}_0(B_\rho(x_0); \bR^N)$ such that
\begin{equation} \label{Ekeland:1}
\mvint_{B_\rho(x_0)} |Du -Dw|^{p_2} dx \leq
C(\bM) \rho^{p_2 \alpha /2},
\end{equation}
and
\begin{equation} \label{Ekeland:2}
\mvint_{B_\rho(x_0)} f(x_m,Dw) dx \leq \mvint_{B_\rho(x_0)}
f(x_m,Dw+D\phi) dx + c\rho^{\alpha /2} \mvint_{B_\rho(x_0)} (1+ |D\phi|^{p_2}) dx,
\end{equation}
for any $\phi \in W^{1,p_2}_0(B_\rho(x_0); \bR^N)$.
\end{proposition}

\begin{proof} Consider the function $g(z) := f(x_m,z)$. Then $g$ satisfies the growth
condition \eqref{growth} with exponent $p_2 = p(x_m)$. Let $v
\in u + W^{1,p_2}_0(B_{\rho};\bR^N)$ be the unique solution of the
Dirichlet problem
\begin{equation*}
     \min\left\{ \int_{B_{\rho}} g(Dw):\ w \in u + W^{1,p_2}_0(B_{\rho};\bR^N)
     \right\}.
\end{equation*}
$v$ exists as $f$ is quasiconvex. Lemma \ref{HI.bdry} with $p\equiv
p_2, q \equiv p_2(1+\delta/4)$ and $a(x) \equiv 1$ provides $\eps
\equiv \eps(\gamma_1, \gamma_2,L)$ and $c \equiv
c(\gamma_1,\gamma_2,L)$ with $0 < \eps < \delta/4$ such that
\begin{multline*}
     \left( \mvint_{B_\rho(x_0)} |Dv|^{p_2(1+\eps)}\dx\right)^{\frac{1}{p_2(
     1+\eps)}}\\
     \le c\left(\mvint_{B_\rho(x_0)} |Dv|^{p_2}\dx\right)^{1/p_2}
     +\left(\mvint_{B_\rho(x_0)} \left(|Du|^{p_2}+1\right)^{1+\delta/4}\dx\right)^{
     \frac{1}{p_2(1+\delta/4)}}.
\end{multline*}
Higher integrability for the function $u$ (Lemma \ref{High.Int.u}) gives us
\begin{equation*}
     \midint_{B_{\rho}} |Du|^{p_2(1+\delta/4)}\dx \le C(\bM).
\end{equation*}
For $v$ we get by the minimality and the growth condition
\begin{equation*}
     \midint_{B_{\rho}} |Dv|^{p_2}\dx \le L^2 \midint_{B_{\rho}} 1+ |Du|^{p_2}\dx,
\end{equation*}
so that together with the estimate before and the higher integrability for $u$ we
have
\begin{equation*}
     \midint_{B_{\rho}} |Dv|^{p_2(1+\eps)}\dx \le  C(\bM),
\end{equation*}
with $0<\eps <\delta /4$ and $C(\bM)$ also depending on
$n,L,\gamma_1,\gamma_2$.
We now estimate the difference
\begin{equation*}
\begin{aligned}
     \midint_{B_{\rho}}( g(Du) - g(Dv))\dx
     & = \midint_{B_{\rho}} \left( f(x_m,Du)-f(x_m,Dv)\right)\dx\\
     & = \midint_{B_{\rho}} \left( f(x,Du)-f(x,Dv)\right)\dx\\
     &\mskip+60mu + \midint_{B_{\rho}} \left(f(x,Dv)-f(x_m,Dv)\right)\dx\\
     &\mskip+60mu + \midint_{B_{\rho}}\left(f(x_m,Du)-f(x,Du)\right)\dx\\
     & = (I) + (II) + (III).
\end{aligned}
\end{equation*}
By the minimality of $u$ we have $(I) \le 0$. For $(III)$ we get by the continuity
of $f$ in the first variable and the higher integrability of $u$
\begin{equation*}
     |(III)| \le c(\delta) \omega(\rho) \midint_{B_{\rho}} (|Du|^{p_2(1+\delta/4)}
     +1)\dx \le C(\bM) \omega(\rho).
\end{equation*}
By the same arguments for the function $v$ we have
\begin{equation*}
     |(II)| \le c(\eps)\omega(\rho) \midint_{B_{\rho}} (|Dv|^{p_2(1+\eps
     )}+1)\dx \le C(\bM) \omega(\rho),
\end{equation*}
so overall, using \eqref{cont.strong} we get
\begin{equation}\label{comp.}
     \midint_{B_{\rho}}[ g(Du) - g(Dv)]\dx \le c \omega(\rho) \le C(\bM)
     \rho^{\alpha},
\end{equation}
with the constant $C(\bM)$ depending additionally on
$n,L,\gamma_1,\gamma_2$. Let $0< \mu< \alpha$, $X := u +
W^{1,p_2}_0(B_{\rho};\bR^N)$ and
\begin{equation*}
d: X \times X \to [0,\infty), \quad
d(z,w) = \frac{1}{C(\bM) \rho^{\mu}} \left(\mvint_{B_{\rho}}
|D(z-w)|^{p_2} \dx\right)^{1/p_2}.
\end{equation*}
On the
complete metric space $(X,d)$ we consider the functional
\begin{equation*}
     \cG:X \to \bR,\ \cG(z) := \midint_{B_{\rho}} g(Dz)\, dx ,
\end{equation*}
which is clearly lower semicontinuous. By $\cG(v) = \min_X \cG$ and \eqref{comp.} we
have
\begin{equation*}
     \cG(u) \le \inf_X \cG + C(\bM) \rho^{\alpha}
\end{equation*}
Therefore the Ekeland variational principle (Lemma \ref{eke_orig})
provides a function $w \in u + W^{1,p_2}_0(B_{\rho};\bR^N)$ with the
properties
\begin{eqnarray*}
     &&\mvint_{B_{\rho}} |Du - Dw|^{p_2}\dx \le C(\bM) \rho^{\mu p_2} \quad
          \mbox{ and } \\
     &&\mvint_{B_{\rho}} g(Dw)\dx \le \mvint_{B_{\rho}} g(Dw + D\varphi)\dx +
          \rho^{\alpha - \mu} \left(\mvint_{B_{\rho}} |D\varphi
          |^{p_2}\dx \right)^{1/p_2},
\end{eqnarray*}
for all $\varphi \in W^{1,p_2}_0(B_{\rho};\bR^N)$. We estimate the second integral of
the second inequality simply applying Bernoulli's inequality:
\begin{equation*}
     \left(\mvint_{B_{\rho}}|D\varphi|^{p_2}\, dx\right)^{1/p_2} \le
     \left(1+ \mvint_{B_{\rho}}|D\varphi|^{p_2}\, dx \right)^{1/p_2}
     \le c(p_2) \mvint_{B_{\rho}}\left(1 + |D\varphi|^{p_2}\right)\, dx.
\end{equation*}
Thus choosing $\mu \equiv \alpha/2$ we obtain the assertion.

\end{proof}

\subsection{Caccioppoli inequality}

\begin{Lem}\label{Caccioppoli}
Let $\bM, \hat{M}>0$. Assume that $u$ is a
local minimizer of the functional $\cF$ with $(|Du|^{p_2})_{x_0,2\rho} \le {\bar
M}$ and  $A\in \Hom(\bR^n;\bR^N)$ with $|A| \leq \hat{M}$. There exist constants
$\rho_0 = \rho_0({\bar M},\hat{M},\alpha)$ and $c_c= c_c(\hat{M})$ such that for
every $\xi \in \bR^N$ and every ball $B(x_0,\rho) \Subset \Omega$ with $\rho \leq \rho_0$
there holds:
\begin{equation}
\begin{aligned}
\mvint_{B_{\rho/2}(x_0)} |V_{p_2}&(Du-A)|^2 \dx\\
     & \leq c_c \left[
          \mvint_{B_\rho(x_0)} \left| V_{p_2}\left(\frac{u-\xi-A(x-x_0)}{\rho}\right)
          \right|^2 dx + \rho^{\alpha /2}\right].
\end{aligned}
\end{equation}
\end{Lem}

\begin{proof}
By Proposition \ref{ekeland} there exists an almost minimizer $w$ of the frozen
functional
$\tilde{\cF}_\Omega[\cdot] \equiv \int_{B_{\rho}(x_0)} f(x_m,\cdot\ )\, dx$,
such that
\begin{equation}\label{u-w}
\mvint_{B_\rho(x_0)} |Du -Dw|^{p_2} dx \leq C(\bar M) \rho^{p_2\alpha/2},
\end{equation}
and $w$ satisfies \eqref{Ekeland:2}, so that by Lemma 3 in \cite{duzaar-g-k} with
$\omega(\rho) = \rho^{\alpha/2}$ we have
\begin{equation} \label{cacc-w}
\begin{aligned}
     \mvint_{B_{\rho/2}(x_0)} |V_{p_2}(&Dw-A)|^2 \dx\\
     &\leq c \mvint_{B_\rho(x_0)} \left| V_{p_2}\left(\frac{w-\xi-
          A(x-x_0)}{\rho}\right) \right|^2 dx + c\rho^{\alpha/2},
\end{aligned}
\end{equation}
where $c=c(\hat{M})$. By Lemma \ref{prop.V} $(ii)$ and H\"older's inequality we have
in the case $p_2 \ge 2$
\begin{eqnarray*}
     &&\mvint_{B_{\rho/2}} |V_{p_2}(Du-Dw)|^2\, dx \\
     &&\mskip+80mu \le \mvint_{B_{\rho/2}} |Du-Dw|^{p_2}\,
           dx + \left(\mvint_{B_{\rho/2}} |Du-Dw|^{p_2}\, dx\right)^{2/p_2} \\
     &&\mskip+80mu \le c\left(\rho^{p_2 \alpha/2} + \rho^{\alpha}\right).
\end{eqnarray*}
In the case $1 < p_2 < 2$ we directly estimate
\begin{equation*}
     \mvint_{B_{\rho/2}(x_0)} |V_{p_2}(Dw-A)|^2\, dx \le \mvint_{B_{\rho/2}(x_0)}
     |Du - Dw|^{p_2}\, dx \le c\rho^{p_2\alpha/2},
\end{equation*}
and thus in every case we have
\begin{eqnarray*}
&&\mvint_{B_{\rho/2}} |V_{p_2}(Du-A)|^2 dx\\
     &&\mskip+70mu \le  c(p_2) \left[ \mvint_{B_{\rho/2}} |V_{p_2}(Du-Dw)|^2 dx +
          \mvint_{B_{\rho/2}} |V_{p_2}(Dw-A)|^2 \, dx \right]\\
     &&\mskip+70mu \le  c(p_2) \left[ \rho^{p_2 \alpha/2} + \rho^{\alpha} +
          \mvint_{B_{\rho/2}} |V_{p_2}(Dw-A)|^2 \dx \right].
\end{eqnarray*}
On the other hand, again by the properties of the function $V_{p_2}$ we estimate
\begin{eqnarray*}
     \lefteqn{\mvint_{B_\rho} 
     \left|V_{p_2}\left(\frac{w-\xi-A(x-x_0)}{\rho}\right)\right|^2
          \, dx}\\
     &&\le c(p_2)\left[ \mvint_{B_{\rho}} \left|V_p\left(
          \frac{u-\xi-A(x-x_0)}{\rho}\right)\right|^2\, dx + \mvint_{B_{\rho}}
          \left|V_{p_2}\left(\frac{w-u}{\rho}\right)\right|^2\, dx\right].
\end{eqnarray*}
Since $u-w \in W^{1,p_2}_0(B_{\rho};\bR^N)$, we apply Poincar\'e's inequality
on the second term of the right hand side, finally obtaining (again using properties
of the function $V_{p_2}$)
\begin{equation*}
     \mvint_{B_\rho} \left|V_{p_2}\left(\frac{u-w}{\rho}\right)\right|^2\, dx
     \le c\mvint_{B_\rho}\left|V_{p_2}(Dw-Du)\right|^2\, dx \le c
     \left( \rho^{p_2\alpha/2} + \rho^{\alpha}\right).
\end{equation*}
Hence from \eqref{u-w} and \eqref{cacc-w} we obtain
\begin{multline*}
\mvint_{B_{\rho/2}(x_0)} |V_{p_2}(Du-A)|^2 dx \\
\le c\left[ \rho^{p_2\alpha/2} + \mvint_{B_\rho(x_0)} \left| V_{p_2}\left(
\frac{u-\xi-A(x-x_0)}{\rho}\right) \right|^2 dx + \rho^{\alpha}\right].
\end{multline*}
The claim follows since $\rho^{p_2\alpha /2} \leq \rho^{\alpha /2}$.
\end{proof}

\subsection{Approximate $\cA$-harmonicity}

\begin{Lem} \label{Approx-A-harm}
Let $\bM, \hat{M}>0$. Assume that $u$ is a
local minimizer of the functional $\cF$ with $(|Du|^{p_2})_{x_0,2\rho} \le {\bar
M}$ and  $A\in \Hom(\bR^n;\bR^N)$ with $|A| \leq \hat{M}$.
There exist a constant $c_e \equiv
c_e(n,N,p_2,L,\bM,\hat{M})$ and a radius $\rho_0=\rho_0(\alpha)$ such that for every ball
$B_{\rho}(x_0) \Subset \Omega$ with $\rho \le \rho_0$ we have
\begin{multline*}
     \Biggl| \midint_{B_{\rho}(x_0)} D^2 f(x_0,A)(Du-A,D\varphi)\dx \Biggr| \\
     \le c_e \left( \Phi^2 + \sqrt{\nu_{\hat{M},x_0}(\Phi)}\Phi + \sqrt{\omega(\rho)} \right)
     \sup\limits_{B_{\rho}(x_0)} |D\varphi|,
\end{multline*}
for all $\varphi \in C^1_0(B_{\rho}(x_0);\bR^N)$.
\end{Lem}

\begin{proof}
First we assume $|D\varphi| \le 1$. Let $0 < s \le 1$. We start by
showing the following inequality:
\begin{equation}\label{Equ.1}
\begin{aligned}
     \midint_{B_{\rho}} &D^2 f(x_0,A)(Du -A,D\varphi)\dx \\
     \ge &  \frac{1}{s}\Biggl[
        \midint_{B_{\rho}} \int_0^s \Big( Df(x_0,Du)-Df(x_0,Du+\tau D\varphi)\Big)
        D\varphi \, d\tau \dx \\
     & + s\midint_{B_{\rho}} \int_0^1
        \Big( D^2 f(x_0,A) -D^2 f(x_0,A+\tau(Du-A))\Big)\Big( Du-A,D\varphi\Big) \, d\tau\dx\\
     & - c\omega(\rho) \Biggr]\\
     = & (A) + (B) + (C).
\end{aligned}
\end{equation}
with $c\equiv c(\bM,p_2,\delta)$.
To see that let us start with the difference
\begin{equation*}
     \midint_{B_{\rho}} (f(x_0,Du) - f(x_0,Du + sD\varphi))\dx.
\end{equation*}
Introducing two additional differences we get
\begin{align*}
     \midint_{B_{\rho}} &(f(x_0,Du) - f(x_0,Du + sD\varphi))\dx
     =\midint_{B_{\rho}} f(x_0,Du) - f(x,Du)\dx\\
    & \qquad \qquad + \midint_{B_{\rho}} f(x,Du) - f(x,Du + sD\varphi)\dx\\
    & \qquad \qquad+ \midint_{B_{\rho}} f(x,Du + sD\varphi) - f(x_0,Du+sD\varphi)\,dx\\
    &\leq \midint_{B_{\rho}} |f(x,Du) - f(x_0,Du)|\dx\\
    & \qquad \qquad + \midint_{B_{\rho}} |f(x,Du + sD\varphi) -
     f(x_0,Du+sD\varphi)|\dx\\
     & = (I) + (II)
\end{align*}
since there holds
$$
\midint_{B_{\rho}} f(x,Du) - f(x,Du + sD\varphi)\dx \le 0,
$$
because of the minimality of the function $u$.
To estimate the first term, we use the continuity condition for $f$ with respect to
the variable $x$ as follows:
\begin{eqnarray*}
     &&\mskip-10mu (I) \le \\
     &&\mskip+20mu \midint_{B_{\rho}} \omega(|x-x_0|) \left[ (1+ |Du|^2)^{p(x)/2} +
     (1+|Du|^2 )^{p(x_0)/2}\right] (1+\log(1+|Du|^2))\dx.
\end{eqnarray*}
By the elementary inequality
\begin{equation*}
     \log(1+|z|^2) \le C(a) |z|^a \quad \mbox{ for all } 0 < a < 1,
\end{equation*}
and the fact that $p_2 \ge p(x)$ for all $x$ we see that
\begin{multline*}
     \left[(1+|Du|^2)^{p(x)/2} + (1+|Du|^2)^{p(x_0)/2}\right](1+\log(1+|Du|^2))
     \le \\
     c(p_2,\delta)(1+|Du|^{p_2(1+\delta/4)}),
\end{multline*}
where $\delta$ is the exponent of Lemma \ref{High.Int.u}. Higher
integrability of $u$ gives us (together with estimate \eqref{Est.p1p2} for the
exponents)
\begin{equation*}
     \midint_{B_{\rho}} |Du|^{p_2(1+\delta/4)}\dx \le \midint_{B_{\rho}} 1+
     |Du|^{p(x)(1+\delta)}\dx \le c \left(\midint_{B_{\rho}}|Du|^{p(x)}+1\dx
     \right)^{1+\delta}.
\end{equation*}
This leads to
\begin{equation*}
     (I) \le c\omega(\rho) (1+\bM)^{1+\delta}.
\end{equation*}
For $(II)$ we follow the same way as for $(I)$, additionally using $|D\varphi| \le
1$: The continuity condition in the variable $x$ together with the estimates for the
exponent $p$ lead to
\begin{equation*}
     (II) \le \omega(\rho) \midint_{B_{\rho}} (1+|Du|^2 + s^2|D\varphi|^2)^{
     p_2/2}(1+\log(1+|Du|^2 + s^2|D\varphi|^2))\dx.
\end{equation*}
Using $|D\varphi| \le 1$ we immediately get
\begin{align*}
     (II) &\le c(p_2)\omega(\rho)\midint_{B_{\rho}} (1+|Du|^{p_2(1+\delta/4)})\dx
     \\
     & \le c(\gamma_2,\delta) \omega(\rho) (1+\bM)^{1+\delta}.
\end{align*}

Finally we conclude
\begin{equation}
     \midint_{B_{\rho}} (f(x_0,Du) - f(x_0,Du + sD\varphi))\dx \le c(\bM,\gamma_2,\delta)
     \omega(\rho),
\end{equation}
and hence
\begin{equation*}
     -\midint_{B_{\rho}} \int_0^s Df(x_0,Du+\tau D\varphi)D\varphi \, d\tau\dx-
     c\omega(\rho)\le 0.
\end{equation*}
Secondly we see that (remark that $\varphi \in C_0^1(B_{\rho},\bR^N)$)
\begin{equation*}
\begin{aligned}
     \midint_{B_{\rho}} \int_0^1 D^2f(x_0,A+ \tau(Du-A))&(Du-A,D\varphi)\, d\tau \,
          dx\\
     &= \midint_{B_{\rho}} (Df(x_0,Du) - Df(x_0,A))D\varphi \dx\\
     &= \midint_{B_{\rho}} Df(x_0,Du)D \varphi \dx.
\end{aligned}
\end{equation*}
It follows then that \eqref{Equ.1} is true.

\medskip

We start now taking a look at the right hand side of \eqref{Equ.1}. In what
follows we will distinguish the sets
\begin{equation*}
     B_{\rho}^+ \equiv B_{\rho} \cap \{x: |Du - A| > 1\} \quad \mbox{ and } \quad
     B_{\rho}^- \equiv B_{\rho} \cap \{x: |Du - A| \le 1\}.
\end{equation*}
Let us remark that by Lemma \ref{prop.V} in the case $|Du-A| > 1$ we have for both $1<p_2<2$
and $p_2 \ge 2$ the estimate
\begin{equation*}
     |Du-A|^{p_2} \le c(\gamma_2,\hat{M})|V_{p_2}(Du-A)|^2 \le c|V_{p_2}(Du)-V_{p_2}(A)|^2.
\end{equation*}
In the case $|Du-A| \le 1$ we obtain for all $p_2 > 1$
\begin{equation*}
     |Du-A|^2 \le c(\gamma_2,\hat{M})|V_{p_2}(Du-A)|^2 \le c|V_{p_2}(Du) - V_{p_2}(A)|^2.
\end{equation*}

We first estimate $|(A)|$.
On the set $B_{\rho}^-$, we put $(A)$ in terms of the second derivative of $f$ by
writing
\begin{equation*}
\begin{aligned}
     \int_0^s (Df(x_0,Du)& - Df(x_0,Du+\tau D\varphi ))D \varphi d\tau\\
     & =\int_0^s \int_0^1 D^2 f(x_0,Du+\sigma \tau D \varphi )(\tau D\varphi,
          D\varphi )\, d\sigma \, d\tau.
\end{aligned}
\end{equation*}
As we are on the set $B_{\rho}^-$, we have $|Du + \sigma \tau D\varphi | \le |Du-A|
+ |A| + |D\varphi| \le 2+\hat{M}$ and therefore $|D^2 f(x_0,Du + \sigma \tau
D\varphi )|\le K(\hat{M})$, so we get
\begin{equation*}
     \frac{1}{s}\Bigl|\int_0^s (Df(x_0,Du) - Df(x_0,Du+\tau D\varphi ))D \varphi
     d\tau\Bigr| \le \frac{s}{2}K(\hat{M}).
\end{equation*}
On the set $B_{\rho}^+$ by \eqref{growth_Df} we get
\begin{equation*}
\begin{aligned}
     \frac{1}{s}\Bigl| \int_0^s (Df(x_0,Du)& - Df(x_0,Du+ \tau D \varphi ))D \varphi
          \, d\tau\Bigr|\\
     &\le \frac{L}{s} \int_0^s \left[ |Df(x_0,Du)| + |Df(x_0,Du+\tau D\varphi
          )|\right]|D \varphi|\, d\tau\\
     &\le c(\gamma_2,L) \left[1+|Du|^{p_2-1}\right].
\end{aligned}
\end{equation*}
Since $|Du-A| > 1$, we have $|Du|^{p_2-1} \le |Du-A|^{p_2-1} + |A|^{
p_2-1} \le |Du-A|^{p_2}+ \hat{M}^{p_2-1} \le (1+\hat{M}^{p_2-1})|Du-A|^{p_2}$.
Therefore the last term of the estimate above can be further estimated by
\begin{equation*}
     c(\gamma_2,\hat{M},L,n) |Du-A|^{p_2} \le c|V_{p_2}(Du)-V_{p_2}(A)|^2.
\end{equation*}
Summing up the arguments before we get
\begin{equation*}
     |(A)| \le c(n,L,\hat{M},\gamma_2) \midint_{B_{\rho}} |V_{p_2}(Du) - V_{p_2}(
     A)|^2\dx + \frac{s}{2}K(\hat{M}) = c \Phi^2  +\frac{s}{2}K(\hat{M}).
\end{equation*}

To estimate $|(B)|$, on the set $B_{\rho}^-$ we use the fact that
$|A+\tau(Du-A)| \le \hat{M}+1$ and by \eqref{cont_D2f} and
\eqref{growth_D2f} we obtain
\begin{equation*}
\begin{aligned}
     |D^2f(x_0,A)& - D^2f(x_0,A+\tau(Du-A))|\\
     & =  \left( |D^2f(x_0,A) - D^2f(x_0,A+\tau (Du-A))|^2\right)^{1/2}\\
     &\le \sqrt{2\sup_{|B|\le \hat{M}+1} |D^2f(B)|\ \nu_{\hat{M},x_0}(|Du-A|)}\\
     & =  \sqrt{2K(\hat{M})} \sqrt{\nu_{\hat{M},x_0}(|Du-A|)}.
\end{aligned}
\end{equation*}
Therefore
\begin{equation*}
\begin{aligned}
     \Bigl|\int_0^1 &\left[D^2f(x_0,A)-D^2f(x_0,A+\tau(Du-A))\right](Du-A,D\varphi)\,
          d\tau\Bigr|\\
     &\le \int_0^1 |D^2f(x_0,A)-D^2f(x_0,A+\tau(Du-A))||Du-A||D\varphi |\, d\tau\\
     &\le \sqrt{2K(\hat{M})}\sqrt{\nu_{\hat{M},x_0}(|Du-A|)}\ |Du-A|\\
     &\le c\sqrt{2K(\hat{M})}\sqrt{\nu_{\hat{M},x_0}(|V_{p_2}(Du) - V_{p_2}(A)|)}
          \ |V_{p_2}(Du) - V_{p_2}(A)|.
\end{aligned}
\end{equation*}
On the set $B_{\rho}^+$ we write
\begin{equation*}
\begin{aligned}
     \Bigl|\int_0^1(D^2&f(x_0,A)-D^2f(x_0,A+\tau (Du-A)))(Du-A,D \varphi )\, d\tau
     \Bigr|\\
     &= \Bigl|D^2f(x_0,A)(Du-A,D\varphi ) + (Df(x_0,A) - Df(x_0,Du))D\varphi \bigr|\\
     & =: (*).
\end{aligned}
\end{equation*}
The first term is estimated by $K(\hat{M}+1)$.
For the second term we use again the growth condition for $Df$ and follow
exactly the same way as above for $|(A)|$ on
the set $B_{\rho}^+$ to get
\begin{equation*}
      (*) \le c |Du-A|^{p_2} \le c|V_{p_2}(Du) - V_{p_2}(A)|^2,
\end{equation*}
with $c \equiv c(n,N,p_2,L,\bM,\hat{M})$. Putting the estimates
together, we deduce
\begin{equation}
     |(B)| \le c \left(\sqrt{\nu_{\hat{M},x_0}(\Phi)} \Phi + \Phi^2\right).
\end{equation}
We choose $s \equiv \sqrt{\omega(\rho)} < 1$ and get
\begin{equation*}
     |(A)| + |(B)| + |(C)| \le c\left(\Phi^2 + \sqrt{\nu_{\hat{M},x_0}(\Phi)}\Phi + \sqrt{
     \omega(2\rho)}\right),
\end{equation*}
with a constant $c \equiv c(n,N,\gamma_2,L,\bM,\hat{M})$.
Altogether we have shown
\begin{equation*}
     \midint_{B_{\rho}} D^2f(x_0,A)(Du-A,D\varphi )\dx \ge -c\Big( \Phi^2 + \sqrt{
     \nu_{\hat{M},x_0}(\Phi)}\Phi + \sqrt{\omega(\rho)}\Big).
\end{equation*}
The estimate for $\midint Df(x_0,A)(Du-A,D\varphi )\dx$ from above
is shown exactly in the same way. This gives the lemma for the case
$|D\varphi| \le 1$ and the general result can be achieved by
rescaling.

\end{proof}

\section{Proof of the result}

\subsection{Excess-improvement lemma} \label{excess-dec}

Let $M > 0$ be fixed. Consider a point $x_0 \in \Omega$ such that
$|(Du)_{x_0, 2\rho}| \le M$ and $\Phi(x_0, 2\rho) \le 1$., where
$\Phi$ is the function defined in \eqref{excess} with $A=(Du)_{x_0,\rho}$,  i.e.
$$
     \Phi(x_0, \rho) = \Phi(x_0,\rho,(Du)_{x_0,\rho})
     = \left( \midint_{B_{\rho}(x_0)} |V_{p_2}(Du) -
     V_{p_2} ((Du)_{x_0,\rho})|^2\dx \right)^{1/2}.
$$
By properties of the function $V_{p_2}$ (Lemma \ref{prop.V}) we obtain
\begin{eqnarray*}
     (|Du|^{p_2})_{x_0,2\rho}
     &\le& 2^{p_2-1} \mvint_{B_{2\rho}} |Du - (Du)_{x_0,2\rho}|^{p_2}
          \, dx + 2^{p_2-1} |(Du)_{x_0,2\rho}|^{p_2}\\
     &\le& 2^{p_2-1}c \mvint_{B_{2\rho}} |V_{p_2}(Du - (Du)_{x_0,2\rho})|^2\, dx\\
     && +
          2^{p_2-1}c\left( \mvint_{B_{2\rho}} |V_{p_2}(Du-(Du)_{x_0,2\rho}|^2\, dx
          \right)^{p_2/2}+ 2^{p_2-1}M^{p_2}\\
     &\le& 2^{p_2-1}c_1(M)\left[ \Phi^{p_2}(x_0,2\rho) + \Phi^2(x_0,2\rho) + M^{p_2}\right]\\
     &\le& 2^{\gamma_2-1}c_1(M) \left[ 2 + M^{\gamma_2} \right] \equiv \bar{M},
\end{eqnarray*}
where $c_1(M)$ is the constant out of Lemma \ref{prop.V}.
Furthermore we remark that by H\"older's inequality we immediately get
\begin{equation}
|(Du)_{x_0,\rho}| \le 2^n (|Du|)_{x_0,2\rho} \le
2^n (|Du|^{p_2})_{x_0,2\rho}^{1/p_2} \le 2^n \bar{M}^{1/p_2} \equiv \hat{M}.
\end{equation}
By the quasiconvexity condition we get that ${\mathcal A} \equiv
D^2f(x_0,(Du)_{x_0,\rho})$ is elliptic in the sense
of Legendre-Hadamard with ellipticity constant $\kappa$ and upper
bound $K$, where
\begin{equation*}
\kappa \equiv  \frac{2}{L}(1+\hat{M}^2)^{(p_2-2)/2} \quad \text{and} \quad
K \equiv K_{\hat{M}} \equiv \sup_{|A|\le \hat{M}} |D^2f(x_0,A)|.
\end{equation*}
Recall now the notation for constants. We will skip their dependence
on $n,N,L$ and $\gamma_1, \gamma_2$ and remark their dependence on $M$ and
$\alpha$. We will not remark the dependence on $\hat{M}$ and $\bar{M}$ since the
two constants are computed out of $M,n,\gamma_1, \gamma_2$. Denote
\begin{eqnarray*}
c_1 &=& c_1(M) \quad \quad \text{from Lemma \ref{prop.V} $(v)$ and $(vi)$ (properties
of function $V$)}\\
c_a &=& c_a(M) \quad \quad \text{from Lemma \ref{A-harm} (estimates for $\cA$--harmonic
function)}\\
c_c &=& c_c(M) \quad \quad \text{from Lemma \ref{Caccioppoli} (Caccioppoli
inequality)}\\
\rho_0 &=& \rho_0(M,\alpha) \quad \mskip-2mu \text{from Lemma \ref{Caccioppoli} (Caccioppoli inequality)}\\
c_e &=& c_e(M) \quad \quad \text{from Lemma \ref{Approx-A-harm}
(approximate $\cA$--harmonicity)}
\end{eqnarray*}

\begin{Lem}\label{excess-imp}
Let $M > 0$ and $\beta \in (\alpha/4,1)$ be fixed.  There exist $\theta =\theta(\beta,M)\in
(0,1/4]$ and $\delta = \delta(\theta) \in (0,1]$ such that
if $x_0 \in \Omega$ is a point such that
\begin{eqnarray*}
     |(Du)_{x_0, 2\rho}| &\le& M, \\
     \sqrt{\nu_{\hat{M},x_0}(\Phi(x_0,\rho))} +\Phi(x_0,\rho) &\le& \delta/2, \\
     2\sqrt{2}c_a c_1 c_e  \sqrt{ \Phi^2(x_0,\rho) + \frac{4}{\delta^2} \omega(\rho)}&\le&
     1,
\end{eqnarray*}
hold for some $\rho \in (0,\rho_0]$, then
\begin{equation*}
\Phi^2(x_0,\theta \rho) \le \theta^{2\beta} \Phi^2(x_0,\rho) + \hat{c}
\rho^{\alpha /2},
\end{equation*}
with $\hat{c} \equiv \hat{c}(M,\theta)$.
\end{Lem}

\begin{proof}
With $x_0$ fixed we will write $\Phi(\rho)$ instead of $\Phi(x_0, \rho)$.
Let $\theta$ be a parameter, which is free at first and will be fixed at the end of
the proof. Set
\begin{equation} \label{eps}
\eps \equiv
\begin{cases}
\theta^{n+4} & \quad \text{if $p_2 \in (1,2)$}, \\
\theta^{n+ p_2 + 2} &\quad \text{if $p_2 \geq 2$},
\end{cases}
\end{equation}
and $\delta \equiv
\delta(\varepsilon) \equiv \delta(\theta)$ be the parameter out of Lemma
\ref{A-approx} (lemma on $\cA$--harmonic approximation).

By Lemma \ref{prop.V} $(v)$ and $(vi)$ and using the fact that $\theta \le 1$ we estimate
\begin{eqnarray} \label{excess-theta-rho}
     \Phi^2(\theta\rho)
     & = & \mvint_{B_{\theta \rho}(x_0)} \left| V_{p_2} (Du) - V_{p_2} ((
          Du)_{x_0, \theta\rho})\right|^2\, dx \nonumber\\
     &\le& c\mvint_{B_{\theta\rho}(x_0)} \left|V_{p_2}(Du - (Du)_{x_0,\rho} - \gamma
          Dh(x_0))\right|^2\, dx \nonumber\\
     &&\mskip+50mu + c \left| V_{p_2}( (Du)_{x_0,\theta
          \rho} - (Du)_{x_0,\rho} - \gamma Dh(x_0))\right|^2,
\end{eqnarray}
with $c$ depending on $n,N,\gamma_1, \gamma_2$ and also on $M$ (the dependence
on $M$ is due to the constant $c_1(M)$).
We denote
\begin{equation*}
     I \equiv \mvint_{B_{\theta\rho}(x_0)} \left| V_{p_2}(Du - (Du)_{x_0,\rho} -
     \gamma Dh(x_0))\right|^2\, dx.
\end{equation*}
To estimate the second expression on the right hand side of \eqref{excess-theta-rho} we
split the domain of integration into the subsets where $|Du - (Du)_{x_0,\rho} -
\gamma Dh(x_0)| \ge 1$ and $|Du - (Du)_{x_0,\rho} - \gamma Dh(x_0)| < 1$.
Applying Lemma \ref{prop.V} $(iv)$ we obtain
\begin{align*}
     \Bigl| (Du)_{x_0,\theta\rho} - (Du)_{x_0,\rho} - \gamma Dh(x_0)\Bigr|
     & \le \mvint_{B_{\theta\rho}(x_0)} \left|Du - (Du)_{x_0,\rho} - \gamma
          Dh(x_0)\right|\, dx\\
     & \le
     \begin{cases}
     c \left( I^{1/p_2} + I^{1/2}\right) & \quad \text{for $p_2\in(1,2)$},\\
     I^{1/2} & \quad \text{for $p_2 \geq 2$}.
     \end{cases}
\end{align*}
Again applying Lemma \ref{prop.V} $(iv)$ we get
\begin{equation*}
\begin{array}{llll}
V^2_{p_2} (I^{1/p_2} + I^{1/2}) & \leq& cI  &\text{for $p_2 \in(1,2)$},\\ \\
V^2_{p_2} (I^{1/2}) &\leq& c(I + I^{p_2/2})  &\text{for $p_2 \geq 2$}.
\end{array}
\end{equation*}
Thus from \eqref{excess-theta-rho} we conclude
\begin{equation} \label{excess-theta-rho2}
     \Phi^2(\theta\rho) \leq
     \begin{cases}
     cI & \quad \text{for $p_2\in(1,2)$}, \\
     c(I + I^{p_2/2} )& \quad \text{for $p_2 \geq 2$},
     \end{cases}
\end{equation}
with $c \equiv c(M)$. In this place we distinguish the cases $p_2\in (1,2)$ and
$p_2 \geq 2$. However, this is only for technical reasons, i.e. due to the difference in
Lemma \ref{prop.V} $(iv)$. We will see later that in fact $I$ is small and therefore
we can skip the term $I^{p_2/2}$ in the above estimate. We proceed then with the estimates for $I$.

\medskip

Set
\begin{eqnarray*}
     \gamma(\rho) &\equiv& c_1 c_e \sqrt{ \Phi^2(\rho) + \frac{4}{\delta^2}
     \omega(\rho)},\\
     w &\equiv& u - (Du)_{x_0,\rho}(x-x_0).
\end{eqnarray*}
By Lemma \ref{prop.V} $(vi)$ and the smallness condition we have
\begin{equation*}
\mvint_{B_\rho(x_0)} |V_{p_2}(Dw)|^2 \dx \le c_1 \Phi^2(\rho) \le
\gamma^2.
\end{equation*}
It follows from Lemma \ref{Approx-A-harm} and the smallness condition that the function
$w$ is approximate $\cA$-harmonic, i.e.
\begin{equation*}
\begin{aligned}
     \Biggl| \mvint_{B_{\rho}(x_0)} D^2 f&(x_0,(Du)_{x_0,\rho})(
     Dw,D\varphi)\dx \Biggr| \\
     &\le c_e \left( \Phi^2(\rho) + \sqrt{\nu_{\hat{M},x_0}(\Phi(\rho))}\Phi(\rho)
          + \sqrt{\omega(\rho)} \right)\sup_{B_{\rho}(x_0)} |D\varphi| \\
     &\le c_e \left(\Phi(\rho) \frac{\delta}{2}+\sqrt{\omega(\rho)}\right)
     \sup_{B_{\rho}(x_0)} |D\varphi| \\
     &\le \gamma \delta \sup_{B_{\rho}(x_0)} |D\varphi|.
\end{aligned}
\end{equation*}
Thus we are in the situation to apply Lemma \ref{A-approx} ($\cA$--harmonic
approximation) providing a function
$h \in W^{1,p_2}(B_{\rho}(x_0);\bR^N)$ which is $D^2 f(x_0,(Du)_{
x_0,\rho})$-harmonic, such that
\begin{equation} \label{est-h}
\mvint_{B_\rho(x_0)} |V_{p_2}(Dh)|^2 dx \le 1 \quad \text{and} \quad
\mvint_{B_\rho(x_0)} \left|V_{p_2}\left(\frac{w-\gamma
h}{\rho}\right)\right|^2 dx \le \gamma^2 \eps.
\end{equation}
Splitting the domain of integration in the first integral of \eqref{est-h}
into the sets $\{ |Dh| \ge 1\}$ and $\{ |Dh| < 1 \}$ and using
Lemma \ref{prop.V} $(iv)$ we obtain the upper bound for the mean value of $|Dh|$:
\begin{equation}\label{Dh}
     \mvint_{B_{\rho}(x_0)} |Dh|\, dx \le 2 \sqrt{2}.
\end{equation}
Lemma \ref{A-harm} provides
\begin{equation} \label{est-h2}
     \rho \sup_{B_{\rho/2}(x_0)} |D^2 h| + \sup_{B_{\rho/2}(x_0)} |Dh| \le
     c_a \mvint_{B_\rho(x_0)} |Dh|\, dx.
\end{equation}
By assumption there holds $(Du)_{x_0,\rho} \le \hat{M}$ and since
by the smallness assumption $2\sqrt{2} c_a \gamma \leq 1$, we conclude with
\begin{equation*}
     |(Du)_{x_0,\rho}| + \gamma |Dh(x_0)| \le \hat{M} + \gamma c_a \mvint_{
     B_\rho(x_0)} |Dh|\, dx \le \hat{M} + 1.
\end{equation*}
Therefore we can apply the Caccioppoli inequality (Lemma \ref{Caccioppoli})
with $\xi \equiv \gamma h(x_0)$ and
$A \equiv (Du)_{x_0,\rho} + \gamma Dh(x_0)$ obtaining
\begin{align} \label{split-I}
I&= \mvint_{B_{\theta\rho}(x_0)} \left|V_{p_2}(Du - (Du)_{x_0,\rho} - \gamma Dh(x_0))
     \right|^2\, dx \nonumber \\
     &\le c_c \left[\mvint_{B_{2\theta\rho}(x_0)} \left|V_{p_2}\left(\frac{w-\gamma h(
          x_0) - \gamma Dh( x_0)(x-x_0)}{2\theta\rho}\right)\right|^2\, dx +
          \rho^{\alpha/2} \right] \nonumber \\
     &\le c \Biggl[ \mvint_{B_{2\theta\rho}(x_0)} \left|V_{p_2}\left(\frac{w-\gamma
          h}{2\theta\rho}\right)\right|^2 \, dx \nonumber \\
     &\qquad + \mvint_{B_{2\theta\rho}(x_0)}\left|V_{p_2}\left(\gamma \frac{h-h(
          x_0) - Dh(x_0)(x-x_0)}{2\theta\rho}\right)\right|^2\, dx +
      \rho^{\alpha/2}
          \Biggr].
\end{align}
The first integral on the right hand side of \eqref{split-I} is estimated via \eqref{est-h} and Lemma \ref{prop.V} $(i)$
\begin{align} \label{I1}
     \mvint_{B_{2\theta\rho}(x_0)} \left|V_{p_2}\left(\frac{w-\gamma h}{2\theta
          \rho}\right)\right|^2\, dx
     &\le(2\theta)^{-n} \mvint_{B_\rho(x_0)}
          \left| V_{p_2}\left(\frac{w-\gamma h}{2 \theta \rho}\right)\right|^2\,
      dx \nonumber \\
     &\le (2\theta)^{-n} \max\left\{ \left(\frac{1}{2\theta}\right)^2,
     \left(\frac{1}{2\theta}\right)^{p_2}\right\} \gamma^2 \eps \nonumber \\
     &= c(n) \theta^2 \gamma^2,
\end{align}
where the last inequality follows by the choice of $\eps$ \eqref{eps}.
To estimate the second integral on the right hand side of \eqref{split-I} observe that
by \eqref{Dh}, \eqref{est-h2}, applying the Taylor theorem to $h$ on
$B_{2\theta\rho}(x_0)$ we obtain
$$
\sup_{B_{2\theta\rho}(x_0)} |h(x) -  h(x_0) - Dh(x_0)(x-x_0)|
\le 4 \sqrt{2} c_a \theta^2 \rho.
$$
Thus taking $\theta$ sufficiently small we have
$$
\left| \gamma \frac{h-h(x_0)-Dh(x_0)(x-x_0)}{2\theta\rho} \right| \leq 2\sqrt{2}
c_a \theta \gamma \leq 1.
$$
Applying Lemma \ref{prop.V} $(iv)$ in both cases $p_2 \in (1,2)$ and $p_2
\geq 2$ we therefore obtain
\begin{align} \label{I2}
 \mvint_{B_{2\theta\rho}(x_0)} \biggl| V_{p_2}&\left( \gamma \frac{h-h(x_0)-Dh(
          x_0)(x-x_0)}{2\theta\rho}\right)\biggr|^2\, dx \nonumber \\
&\le \mvint_{B_{2\theta\rho}(x_0)} \left| \gamma \frac{h-h(x_0) -Dh(x_0)(
          x-x_0)}{2\theta\rho}\right|^2\, dx \nonumber \\
& \le 8 c_a^2 \theta^2 \gamma^2.
\end{align}
Alltogether, by \eqref{split-I}, \eqref{I1} and \eqref{I2} we end up with
\begin{equation*}
     I \le c\left[c\theta^2\gamma^2 + 8c_a^2\theta^2\gamma^2 + \rho^{\alpha/2} \right].
\end{equation*}
Taking $\theta$ and $\rho$ smaller if necessary we can have $I < 1$.
Therefore, since
$$
\gamma^2 \equiv c_1^2 c_e^2 \left( \Phi^2(\rho) + \frac{4}{\delta^2}
     \rho^\alpha \right),
$$
it follows from \eqref{excess-theta-rho2} that
\begin{equation*}
     \Phi^2(\theta\rho) \le \tilde{c}\theta^2 \Phi^2(\rho) + \hat{c}
     \rho^{\alpha/2},
\end{equation*}
with constants $\tilde{c} \equiv
\tilde{c}(M)$ and $\hat{c} \equiv
\hat{c}(M,\theta,\delta)$ (obviously the constants depend also on the structural
constants $n,N,L,\gamma_1,\gamma_2$).
We now fix $\theta$ sufficiently small, so that
\begin{equation*}
\tilde{c}\theta^2 \le \theta^{2 \beta}.
\end{equation*}
Choice of $\theta$ fixes $\epsilon$ and $\delta$ and the claim follows.

\end{proof}

\subsection{H\"older continuity of $Du$ and a regular set}

Let $M, \theta, \delta$ be fixed. If we assume that $\eta >0$ is such that
\begin{eqnarray*}
     \sqrt{\nu_{\hat{M},x_o}(\eta)} +\eta &\le& \delta/2, \\
     2\sqrt{2}c_1 c_a c_e  \sqrt{ \eta^2 + \frac{4}{\delta^2} \omega(\rho)}&\le&
     1, 
\end{eqnarray*}
and moreover it satisfies some additional technical smallness conditions, and
$\rho_1$ is also sufficiently small, then,  by a standard iteration technique
one obtains 
\begin{equation}
\Phi^2(\theta^j \rho) \le \theta^{2\beta j} \Phi^2(\rho) +
c(M,\alpha, \theta, \delta) (\theta^j \rho)^{\alpha /2},
\end{equation}
for $j = 1,2, \ldots$ provided 
\begin{eqnarray*}
2 \rho &\le& \rho_1, \\
|(Du)_{x_0,2\rho}| &\le& \frac{M}{2},\\
\Phi(\rho) &\le& \frac{\eta}{\sqrt{2}},\\
\Phi(2\rho) &\le& \frac{\eta}{\sqrt{2}}.
\end{eqnarray*}

We define
\begin{eqnarray*}
     &&\Sigma_1 \equiv \left\{ x_0 \in \Omega:\ \limsup_{\rho \downarrow 0} 
     |(Du)_{x_0,\rho}| < + \infty \right\},\\
     &&\Sigma_2 \equiv \left\{ x_0 \in \Omega:\ \limsup_{\rho \downarrow 0} \mvint_{
          B_\rho(x_0)} |Du - (Du)_{x_0,\rho}|^p\, dx = 0 \right\}.
\end{eqnarray*}
The set $\Sigma_1 \cap \Sigma_2$ is a set of full Lebesgue measure,
i.e. ${\mathcal L}^n\left(\Omega \setminus (\Sigma_1 \cap
\Sigma_2)\right) = 0$. The assumptions of Lemma
\ref{excess-imp} as well as the above conditions needed for iteration are satisfied in points $x_0 \in \Sigma_1 \cap
\Sigma_2$.
By standard interpolation we obtain an excess--decay estimate
\begin{equation*}
     \Phi(_0,r,(Du)_r) \le c\left[ \left(\frac{r}{\rho}\right)^{\beta} \Phi(
     x_0,\rho,(Du)_\rho) + \rho^{\alpha/4}\right],
\end{equation*}
for all $0 \le r \le \rho$ and with $\beta\in (\frac{\alpha}{4},1)$.
The regularity result then follows from the fact that this excess--decay estimate implies
\begin{eqnarray*}
     \mvint_{B_r(x)} |V_{p_2}(Du) - (V_{p_2}(Du))_{x,r}|^2\, dy
     &\le& \mvint_{B_r(x)} |V_{p_2}(Du) - V_{p_2}((Du)_{x,r})|^2\, dy\\
     &\le& c \left[ \left(\frac{r}{\rho}\right)^{\beta} \Phi(x_0,\rho,(Du
          )_\rho) + \rho^{\alpha/4} \right],
\end{eqnarray*}
for any $x$ in the neighbourhood of $x_0$. From this estimate we conclude, by
Campanato's characterization of H\"older continuous functions (see
\cite{campanato}, \cite{campanato2}), that $V_{p_2}(Du)$ is H\"older continuous with
the exponent $\frac{\alpha}{4}$ in a neighbourhood of $x_0$.

\bigskip

In order to pass over from H\"older continuity of the function $V_{p_2}(Du)$ to the
function $Du$ itself, we use the following 
\begin{Lem}\label{Stetigkeit}
Let $p > 1$ and $w: B \to \bR^N$ a function such that the function
$V_p \circ w: B \to \bR^N$ is H\"older continuous with an exponent $\alpha$. 
Then $w$ is H\"older continuous with the exponent $\beta := \min \left\{1, 2/p\right\} \alpha$.
\end{Lem}

From the above lemma we obtain that $Du$ is locally H\"older continuous with the
exponent
$$
\min \left\{1, 2/\gamma_2\right\} \frac{\alpha}{4}.
$$


\bibliographystyle{plain}
\label{bibliography}
\makeatletter
\addcontentsline{toc}{chapter}{\bibname}
\makeatother
\bibliography{Aharm-px}

\begin{thebibliography}{10}

\bibitem{acerbi-f}
E.~Acerbi and N.~Fusco.
\newblock A regularity theorem for minimizers of quasiconvex integrals.
\newblock {\em Arch. Ration. Mech. Anal.}, 99:261--281, 1987.

\bibitem{acerbi-m1}
E.~Acerbi and G.~Mingione.
\newblock {Regularity results for a class of functionals with nonstandard
  growth.}
\newblock {\em Arch. Ration. Mech. Anal.}, 156(2):121--140, 2001.

\bibitem{acerbi-m2}
E.~Acerbi and G.~Mingione.
\newblock Regularity results for a class of quasiconvex functionals with
  nonstandard growth.
\newblock {\em Ann. Scuola Norm. Sup. Pisa Cl. Sci. IV}, 30(2):311--339, 2001.

\bibitem{acerbi-m3}
E.~Acerbi and G.~Mingione.
\newblock {Regularity results for stationary electro-rheological fluids.}
\newblock {\em Arch. Ration. Mech. Anal.}, 164(3):213--259, 2002.

\bibitem{acerbi-m4}
E.~Acerbi and G.~Mingione.
\newblock {Gradient estimates for the $p(x)$- Laplacean system.}
\newblock {\em J. Reine Angew. Math.}, 584:117--148, 2005.

\bibitem{campanato}
S.~Campanato.
\newblock {Propriet\`a di una famiglia di spazi funzionali}.
\newblock {\em Ann. Scuola Norm. Sup. Pisa}, 18(3):137--160, 1964.

\bibitem{campanato2}
S.~Campanato.
\newblock {Equazioni ellitichi del $II^e$ ordine e spazi
  $\mathcal{L}^{2,\lambda}$}.
\newblock {\em Ann. Mat. Pura Appl.}, 69:321--382, 1965.

\bibitem{carozza-f-m}
M.~Carozza, N.~Fusco, and G.~Mingione.
\newblock {Partial regularity of minimizers of quasiconvex integrals with
  subquadratic growth.}
\newblock {\em Ann. Mat. Pura Appl.}, 175(1):141--164, 1998.

\bibitem{duzaar-g-k}
F.~Duzaar, J.~Grotowski, and M.~Kronz.
\newblock {Regularity of almost minimizers of quasi--convex variational
  integrals with subquadratic growth.}
\newblock {\em Ann. Mat. Pura Appl.}, 184(4):421--448, 2005.

\bibitem{duzaar-s}
F.~Duzaar and K.~Steffen.
\newblock {Optimal interior and boundary regularity for almost minimizers to
  elliptic variational integrals.}
\newblock {\em J. Reine Angew. Math.}, 546:73--138, 2002.

\bibitem{ekeland}
I.~Ekeland.
\newblock {Nonconvex minimization problems.}
\newblock {\em Bull. Am. Math. Soc., New Ser.}, 1:443--474, 1979.

\bibitem{giaquinta}
M.~Giaquinta.
\newblock {\em Multiple Integrals in the Calculus of Variations and Nonlinear
  Elliptic Systems}.
\newblock Princeton University Press, 1983.

\bibitem{jens}
J.~Habermann.
\newblock {Regularity results for functionals and Calder\'on-Zygmund estimates
  for systems of higher order with $p(x)$ growth}.
\newblock Ph.D Thesis, 2006.

\bibitem{zhikov}
V.V. Zhikov.
\newblock {On {L}avrentiev's phenomenon}.
\newblock {\em {Russian J. Math. Phys.}}, 3:249--269, 1995.

\bibitem{zhikov2}
V.V. Zhikov.
\newblock {On some variational problems}.
\newblock {\em {Russian J. Math. Phys.}}, 5:105--116, 1997.

\end{thebibliography}
\nocite{*}

\end{document}